\documentclass[10pt]{article}

\usepackage[francais]{babel}
\usepackage[utf8]{inputenc} 
\usepackage[T1]{fontenc}

\usepackage{amsmath}
\usepackage{amsfonts}
\usepackage{amssymb}

\usepackage{lmodern}
\usepackage[a4paper]{geometry}
\usepackage{graphicx}
\usepackage{amsmath,amscd}
\usepackage{amssymb}
\usepackage{mathtools}
\usepackage{stmaryrd}
\usepackage{calrsfs}
\usepackage{pdfsync}
\usepackage[final]{pdfpages}
\usepackage{enumitem}
\usepackage[francais]{babel} 
\usepackage{esvect}
\usepackage{textcomp}
\usepackage{calc}
\usepackage{verbatim}
\usepackage[colorlinks=true,urlcolor=blue]{hyperref}
\usepackage{fancyhdr}
\usepackage{array}
\usepackage{appendix}
\usepackage{rotating}
\usepackage{multirow}
\usepackage{longtable}
\usepackage{caption} 
\usepackage[all]{xy}
\usepackage{romanbar}

\def \C{{\mathbb C}}
\def \D{{\mathbb D}}
\def \E{{\mathbb E}}

\def \N{{\mathbb N}}

\def \P{{\mathbb P}}

\def \R{{\mathbb R}}

\def \CR{{\mathcal C}}

\def \FR{{\mathcal F}}

\def \v#1{\vspace{#1mm}}
\def \n{{\noindent}}
\def \qeq{{\quad{\rm et}\quad}}
\def \o#1{\overline{#1}}

\def \la{\langle}
\def \ra{\rangle}

\def \vareps{\varepsilon}

\def \norm#1#2{\left\| #1 \right\|_{#2}}

\def \1{{\bf 1}}

\def \begep{\begin{enumerate}[label=\arabic*)]}
\def \bege{\begin{enumerate}}
\def \ende{\end{enumerate}}

\def \var{{\rm var}}
\def \cov{{\rm cov}}

\def \niR{\v5\n{\it Remarque -- }}

\def \niE{\v5\n{\it Exemple -- }}

\def \niP{\v2\n{\it Preuve -- }}

\def \ghyp{g^{\bf h}}
\def \dhyp{d^{\bf h}}
\def \omegahyp{{\omega^{\bf h}}}
\def \nablahyp{{\nabla^{\bf h}}}
\def \starhypb{{*_{\bf h}}}

\def \Deltahyp{{\Delta^{\bf h}}}
\def \hyp {{\bf h}}

\def \cqfd {{\hfill$\footnotesize{\square}$}}

\title{Moyenner sur les cercles le champ libre gaussien\\ dans le disque de Poincaré\\$\quad$ \\
Averaging over the circles the gaussian free field\\ in the Poincaré disk}

\author{Jean-Marc Derrien\footnote{Univ Brest, CNRS UMR 6205, Laboratoire de Mathématiques de Bretagne Atlantique, France. Courriel: jean-marc.derrien@univ-brest.fr}}

\date{9 Janvier 2025}

\begin{document}

\maketitle

\renewcommand{\abstractname}{Résumé}
\begin{abstract}
Le champ libre gaussien sur le disque unité $\D$ peut être appréhendé comme une version bi-dimensionnelle du pont brownien sur l'intervalle $[0,1]$.
Il est intrinsèquement associé à l'espace de Sobolev $H_0^1 (\D)$.
Pour définir ce dernier, on peut choisir n'importe quelle métrique conformément équivalente à la métrique euclidienne sur $\D$.
Cette note est une introduction au
champ libre gaussien sur le disque unité
dont
le but est de mettre en évidence quelques commodités qu'offre la géométrie hyperbolique
sur $\D$ pour décrire les premières propriétés de cet objet probabiliste.
\end{abstract}

\renewcommand{\abstractname}{Abstract}
\begin{abstract}
The gaussian free field on the unit disk $\D$ can be seen as a two-dimensional version of the Brownian bridge on the interval $[0,1]$.
It is intrinsically associated with the Sobolev space $H_0^1 (\D)$.
To define the latter, we can choose any metric conformally equivalent to the Euclidean metric on $\D$.
This note is an introduction to the gaussian free field on the unit disk
whose aim is to highlight some of the conveniences offered by hyperbolic geometry
on $\D$ to describe the first properties of this probabilistic object.
\end{abstract}




\section{Introduction}

Le champ libre gaussien est un objet probabiliste qui peut être défini sur un espace géométrique
discret ou continu
et qui joue un rôle important dans de nombreux domaines.
C'est par exemple une brique essentielle dans la construction de volumes et de métriques 
aléatoires naturels sur l'espace géométrique considéré \cite{Gwynne-Miller}.
Du point de vue de la physique,
c'est aussi un candidat naturel comme objet limite dans une grande variété de situations
en mécanique statistique
et il est aussi présent en théorie quantique des champs \cite{Garban}.
On peut encore mentionner que sa version discrète est un outil puissant pour étudier les temps de recouvrement des marches aléatoires sur les graphes finis (\cite{Ding-Lee-Peres}).

Dans cette note, on s'intéresse plus particulièrement au champ libre gaussien sur le disque unité
$\D=\{z \in \C\, |\ |z|<1\}$ du plan complexe avec l'article \cite{Sheffield} de Sheffield comme référence de base.

Notons $C_c^\infty (\D)$
l'espace vectoriel\footnote{tous les espace vectoriels considérés dans cette note sont des espaces vectoriels sur $\R$} des fonctions définies et régulières\footnote{régulière signifie dans cette note de classe $C^\infty$} sur $\D$ à support compact.
Le champ libre gaussien sur $\D$ est étroitement associé à la structure préhilbertienne
sur $C_c^\infty (\D)$ donnée par le produit scalaire:
$$
\la \varphi,\psi \ra=\int_\D \nabla \varphi (x,y) \cdot \nabla \psi (x,y)\, dxdy
\, .
$$
C'est un objet probabiliste qui vise à palier l'absence de variable aléatoire
qui serait à valeurs dans $C_c^\infty (\D)$, voire dans son complété $H_0^1 (\D)$, 
et dont la projection orthogonale sur n'importe quel sous-espace vectoriel de dimension finie serait de loi normale standard sur ce sous-espace, comme il se doit si l'on ne veut privilégier aucune direction particulière.
Le champ libre gaussien sur $\D$ se définit donc sur la base de propriétés que possèderait une telle variable aléatoire.
Bien que l'on ne puisse parler de sa valeur en un point donné du disque,
ses moyennes sur les cercles de $\D$ prennent un sens naturel.
De plus, la famille de ces moyennes est un processus gaussien dont il existe une version continue.
Il sert de pis-aller pour rendre compte des {\it valeurs} du champ libre gaussien en les points du disque unité. Il permet notamment d'envisager une notion de covariance entre deux telles {\it valeurs}. De plus, lorsqu'il est restreint aux moyennes sur les cercles ayant un même centre fixé, ce processus gaussien réalise un mouvement brownien réel standard dont le temps est paramétré via les rayons des cercles considérés.

Munir le disque unité de la métrique hyperbolique offre quelques avantages dans ces questions.
En effet, par transitivité du groupe des isométries du disque de Poincaré\footnote{c'est-à-dire du disque unité muni de la métrique hyperbolique}, aucun point de $\D$ ne joue de rôle privilégié
et l'on peut au besoin se baser en 0.
De plus,
le produit scalaire sur $C_c^\infty (\D)$ introduit ci-dessus est invariant sous l'action de ces isométries
et
la fonction de Green sur $\D$ utilisée dans l'expression des moyennes considérées
s'exprime simplement à l'aide de la distance hyperbolique.
Enfin, puisque la métrique hyperbolique permet de considérer des cercles dans $\D$ de rayon arbitrairement grand, c'est l'intégralité du mouvement brownien standard dont il est question dans le paragraphe précédent
que l'on réalise en considérant l'ensemble des cercles hyperboliques centrés en un point quelconque de $\D$.

\v2
Voici comment s'organise la suite de cette note.
Dans la section 2, on rappelle quelques notions élémentaires de géométrie différentielle et de géométrie riemannienne dans le plan complexe euclidien qui serviront à établir dans la section 3 des résultats classiques
concernant les moyennes de fonctions sur les cercles, les fonctions harmoniques et l'équation de Poisson.
L'appareillage technique des sections 2 et 3 peut sembler superflu mais il a pour but d'anticiper les besoins des sections 4 et 5, similaires aux deux précédentes mais où le disque $\D$ est cette fois muni de la métrique hyperbolique. La section 4 comprend également un petit développement sur différentes propriétés des isométries dans ce contexte hyperbolique.
Dans la section 6, on considère l'espace de Sobolev $H_0^1 (\D)$ pour lequel on détaille notamment une démonstration de l'invariance de son produit scalaire sous l'action des isométries du disque de Poincaré; puis on s'intéresse, dans la section 7, aux moyennes des éléments de $H_0^1 (\D)$ sur les cercles de $\D$ pour la métrique hyperbolique.
Enfin, dans la section 8, le champ libre gaussien sur $\D$ est défini après quelques rappels sur la notion de variable gaussienne standard à valeurs dans un espace euclidien; puis on discute, dans la dernière section, différentes propriétés du processus constitué de ses moyennes sur les cercles de $\D$ pour la métrique hyperbolique.

\section{Un peu de géométrie dans le plan complexe}

\subsection{Un peu de géométrie différentielle dans le plan complexe}

On suppose le plan complexe $\C$ muni de sa structure de variété différentiable habituelle caractérisée par l'unique carte:
$$
\begin{array}{rccc}
\varphi : & \C & \longrightarrow & \R^2\\
& z=x+i y & \longmapsto & (x(z),y(z))=(x,y)
\end{array}
$$

Une fonction $f$ régulière sur un ouvert $U$ de $\C$
est une fonction $f:U\to\R$ telle que la fonction $f\circ \varphi^{-1}$ soit de classe $C^\infty$
sur l'ouvert $\varphi(U)$ de $\R^2$.

Fixons un point $z$ de $\C$.
L'ensemble des fonctions régulières sur un voisinage de $z$
est naturellement muni d'une structure d'algèbre commutative sur $\R$.
On appelle dérivation en $z$ toute forme linéaire sur cet espace qui vérifie de plus
l'identité de Leibniz: $v(f g)=f(z) v(g)+v(f) g(z)$.
Le plan tangent $T_z \C$ à $\C$ en $z$
est l'espace vectoriel des dérivations en $z$.
\'Etant donnée une fonction $f$ régulière au voisinage de $z$,
on obtient une forme linéaire $df_{|z}$ sur $T_z \C$
en posant:
$$
df_{|z} (v)=v(f)\quad\quad\mbox{pour tout $v$ dans $T_z \C$}\, .
$$
La formule de Taylor à l'ordre 2
assure que
les dérivations partielles
${\partial_x}_{|z}$ et ${\partial_y}_{|z}$ définies par:
$$
{\partial_x}_{|z} f=\frac{\partial f \circ \varphi^{-1}}{\partial x} (\varphi (z))
\qeq
{\partial_y}_{|z} f=\frac{\partial f \circ \varphi^{-1}}{\partial y} (\varphi (z))
$$
constituent une base de $T_z \C$
dont la base duale est donnée par ${dx}_{|z}$ et ${dy}_{|z}$.

Comme il est d'usage, on suppose dans la suite le plan complexe $\C$ orienté par la forme volume $\omega=dx\wedge dy$.

\subsection{Un peu de géométrie dans le plan complexe euclidien}

La métrique euclidienne $g$ sur $\C$
est définie par:
$$
g_z (u_x {\partial_x}_{|z}+u_y {\partial_y}_{|z},v_x {\partial_x}_{|z}+v_y {\partial_y}_{|z})=u_x v_x+u_y v_y
$$
pour tous vecteurs $u=u_x {\partial_x}_{|z}+u_y {\partial_y}_{|z}$ et $v=v_x {\partial_x}_{|z}+v_y {\partial_y}_{|z}$
de $T_z \C$, et pour tout $z$ dans $\C$.
En particulier, $(\partial_x,\partial_y)$ est un champ de bases orthonormées pour cette métrique
et $\omega=dx\wedge dy$ est la forme volume canonique du plan complexe euclidien orienté.

Pour tout $z$ dans $\C$, on définit un isomorphisme $u\mapsto u^*$ de $T_z \C$
dans son dual $T_z^* \C$
en posant:
$$
u^* (v)=g_z (u,v)\quad\quad \mbox{pour tout $v$ dans $T_z \C$}\, .
$$
\'Etant donnée une fonction $f$ régulière sur un ouvert de $\C$,
l'isomorphisme précédent permet de définir, point par point, un champ de vecteurs $\nabla f$ associé à la forme différentielle $df$: $df=(\nabla f)^*$.
Le champ de vecteurs $\nabla f$ ainsi obtenu est le gradient (euclidien) de $f$;
il est régulier sur l'ouvert de $\C$ où $f$ est régulière.

Pour tout $z$ dans $\C$,
on définit un automorphisme
$*:T_z^* \C\to T_z^* \C$
par la relation:
$$
u^* \wedge * \beta = \beta (u) \, \omega_{|z}\quad\quad \mbox{pour tout $u$ dans $T_z \C$ et tout $\beta$ dans $T_z^* \C$}\, .
$$
Cet automorphisme coïncide avec la restriction à $T_z^* \C$ de l'opérateur étoile de Hodge.
Il permet d'associer point par point à toute forme différentielle de degré 1 régulière sur un ouvert de $\C$
une nouvelle forme différentielle de degré 1 également régulière sur le même ouvert de $\C$.
On a en particulier: $*dx=dy$ et $*dy=-dx$.
On définit le laplacien (euclidien) $\Delta f$ d'une fonction $f$ régulière sur un ouvert de $\C$ par la relation:
$$
d*\!df=\Delta f\, \omega
\, .
$$

\niE
Considérons
la fonction $G_0=-\frac{1}{2} \ln(x^2+y^2)$ qui jouera un rôle important dans la suite.
Elle est régulière sur $\C \setminus \{0\}$ et l'on a:
$$
dG_0=-\frac{xdx+ydy}{x^2+y^2}
\, ,
\quad\mbox{et donc aussi}\, ,\quad
*dG_0=-\frac{xdy-ydx}{x^2+y^2}\, .
$$
En différenciant une nouvelle fois, il vient:
$$
d*\!dG_0=
-
\left(
\partial_x \left[ \frac{x}{x^2+y^2}\right]
+
\partial_y \left[ \frac{y}{x^2+y^2}\right]
\right)
dx\wedge dy
=0
\, ;
$$
d'où l'on déduit que $G_0$ est harmonique sur $\C \setminus \{0\}$ au sens où
son laplacien s'annule
sur $\C \setminus \{0\}$.

\section{Moyenner sur les cercles les fonctions régulières dans le plan complexe euclidien}

Dans cette section, le plan complexe orienté $\C$ est supposé muni de la métrique euclidienne $g$.
On note encore $\omega=dx\wedge dy$ sa forme volume canonique.

\v2
Fixons un nombre réel $r$ strictement positif et considérons le cercle $\CR (0,r)$ centré en $0$ de rayon $r$,
muni de la structure de variété riemannienne orientée induite par celle de $\C$.

La moyenne sur $\CR (0,r)$ d'une fonction $f$ mesurable sur ce cercle
est donnée par
$$
f_r (0)
:=
\frac{\int_{\CR (0,r)} f\, \omega_{\CR (0,r)}}{\int_{\CR (0,r)} \omega_{\CR (0,r)}}
=
\frac{1}{2\pi r}
\int_a^{b} f(\gamma (t))\,
\sqrt{g_{\gamma (t)} \left(\dot \gamma (t),\dot \gamma (t)\right)}\, dt
\, ,
$$
où $\omega_{\CR (0,r)}$ désigne la forme volume canonique sur $\CR (0,r)$
et $\gamma:]a,b[\to \C$ n'importe quelle paramétrisation régulière de $\CR (0,r)\setminus\{r\}$\footnote{Par exemple,
avec la paramétrisation: $\gamma (t)=re^{i t}$, $t\in ]0,2\pi[$, on obtient l'expression:
$$
f_r (0)
=
\frac{1}{2\pi}
\int_{0}^{2\pi} f (re^{it})\, dt
\, .
$$}.
En posant, pour une telle paramétrisation $\gamma$,
$\gamma_1 (t):=x(\gamma (t))$
et
$\gamma_2 (t):=y(\gamma (t))$,
on obtient:
$$
f_r (0)
=
\frac{1}{2\pi r}
\int_a^{b} f(\gamma (t))\,
\sqrt{\gamma_1' (t)^2+\gamma_2' (t)^2}
\, dt\, .
$$
Or, pour tout $t$ dans $]a,b[$, on a l'égalité: $\gamma_1 (t)^2+\gamma_2 (t)^2=r^2$, et donc aussi:
$\gamma_1 (t) \gamma_1' (t)+\gamma_2 (t) \gamma_2' (t)=0$.
Il en résulte que si $\gamma$ parcourt le cercle positivement alors il vient:
$$
f_r (0)
=
\frac{1}{2\pi r}
\int_a^{b} f(\gamma (t))\,
\frac{\gamma_1 (t)\gamma_2' (t)-\gamma_1' (t)\gamma_2 (t)}{\sqrt{\gamma_1 (t)^2+\gamma_2 (t)^2}}
\, dt
=
\frac{1}{2\pi}
\int_{\CR (0,r)}
f \, \frac{xdy-ydx}{x^2+y^2}\, .
$$
Ainsi, la moyenne de $f$ sur $\CR (0,r)$ peut s'écrire
de la façon suivante:
$$
f_r (0)
=
-
\frac{1}{2\pi}
\int_{\CR (0,r)}
f \, *\!dG_{0}\, ,
$$
où la fonction $G_0$ a été introduite à la section précédente. 

\v2
Considérons à présent deux nombres réels
$r$ et $R$ avec $0<r<R$. Notons $\o{D} (0,r:R)$
la couronne fermée centrée en $0$ de rayons $r$ et $R$,
et $\partial \o{D} (0,r:R)$ son bord canoniquement orienté.

D'après ce qui précède, pour toute fonction $f$ mesurable sur les cercles $\CR (0,r)$ et $\CR (0,R)$, on a:
$$
f_r (0)
-
f_R (0)
=
\frac{1}{2\pi}
\int_{\partial \o{D} (0,r:R)}
f \, *\!dG_{0}\, .
$$

\v2
On suppose à présent la fonction $f$ régulière
sur un ouvert contenant la couronne $\o{D} (0,r:R)$.

\v1
Grâce à la formule de Green-Riemann, on déduit de l'égalité précédente que
$$
f_r (0)
-
f_R (0)
=
\frac{1}{2\pi}
\int_{D (0,r:R)}
d\left(f \, *\!dG_0 \right)
\, .
$$
Il vient ainsi
\begin{eqnarray*}
f_r (0)
-
f_R (0)
&=&
\frac{1}{2\pi}
\int_{\o{D} (0,r:R)}
\left(df \wedge *dG_0 + f\, d*\!dG_0\right)
\\
&=&
\frac{1}{2\pi}
\int_{\o{D} (0,r:R)}
df \wedge *dG_0 +
\frac{1}{2\pi}
\int_{\o{D} (0,r:R)}
f\, \Delta G_0\, \omega
\\
&=&
\frac{1}{2\pi}
\int_{\o{D} (0,r:R)}
df \wedge *dG_0\, ,
\end{eqnarray*}
où l'on a utilisé que la fonction $G_0$ est harmonique sur $\C\setminus\{0\}$.
Par définition de l'opérateur étoile de Hodge, on a donc:
$$
f_r (0)
-
f_R (0)
=
\frac{1}{2\pi}
\int_{\o{D} (0,r:R)}
g (\nabla f,\nabla G_0)\, \omega\, .
$$

Par les mêmes arguments,
il vient aussi:
\begin{eqnarray*}
& &
\int_{\o{D} (0,r:R)}
g (\nabla f,\nabla G_0)\, \omega
=
\int_{\o{D} (0,r:R)}
dG_0 \wedge *df
\\
&=&
\int_{\o{D} (0,r:R)}
d(G_0 \, *\!df)
-
\int_{\o{D} (0,r:R)} G_0 \, d*\!df
\\
&=&
\int_{\partial \o{D} (0,r:R)} G_0 \, *\!df
-
\int_{\o{D} (0,r:R)} G_0 \,\Delta f\, \omega
\\
&=&
\ln (r) \int_{\CR (0,r)} *df
-
\ln (R) \int_{\CR (0,R)} *df
-\int_{\o{D} (0,r:R)} G_0 \,\Delta f\, \omega\, .
\end{eqnarray*}

\newpage
En définitive, on a donc établi le résultat suivant où sont reprises les notations précédentes.

\v1
\n{\bf Propriété 1 --}
Soient $r$ et $R$ deux nombres réels avec $r<R$.

Si $f$ est régulière sur un ouvert contenant la couronne $\o{D} (0,r:R)$
alors on a:
\begin{eqnarray*}
f_r (0)
-
f_R (0)
&=&
\frac{1}{2\pi}
\int_{\o{D} (0,r:R)}
g (\nabla f,\nabla G_0)\, \omega
\\
&=&
\frac{\ln (r)}{2\pi} \int_{\CR (0,r)} *df
-
\frac{\ln (R)}{2\pi} \int_{\CR (0,R)} *df
-
\frac{1}{2\pi}
\int_{\o{D} (0,r:R)} G_0 \, \Delta f\, \omega
\, .
\end{eqnarray*}

\v5
On en déduit la propriété suivante.

\v1
\n{\bf Propriété 2 --}
Si $f$ est une fonction régulière et harmonique sur le disque ouvert $D (0,R)$ centré en 0 et de rayon $R>0$ et si $f$ est continue sur l'adhérence de ce disque
alors on a:
$$
f (0)=f_R (0)\, .
$$

\niP
Grâce à la formule de Green-Riemann, il vient, pour tout $r$ dans $]0,R[$,
$$
\int_{\CR (0,r)} *df
=
\int_{D (0,r)} d*\!df=\int_{D (0,r)} \Delta f\, \omega=0\, .
$$
Il résulte donc de la propriété 1 que: $f_{R'} (0)=f_r (0)$ pour tous $0<r<R'<R$. On conclut en faisant tendre $r$ vers 0 et $R'$ vers $R$.\cqfd

\v5
On déduit encore de la propriété 1 le résultat suivant qui donne une formule d'inversion (en 0) du laplacien pour les fonctions régulières sur $\C$ à support compact.
 
\v1
\n{\bf Propriété 3 --}
Si $f$ est une fonction régulière sur $\C$ à support compact alors on a:
$$
f (0)=-\frac{1}{2\pi}\int_{\C} G_0 \,\Delta f\, \omega\, .
$$

\niP
Comme $f$ est à support compact,
pour $R$ suffisamment grand
et pour tout $r$ dans $]0,R[$, il vient:
$$
f_R (0)=0\, ,\quad
\int_{\CR (0,R)} *df=0
\qeq
\int_{D (0,r:R)} G_0 \, \Delta f \, d\lambda
=
\int_{\C\setminus D(0,r)} G_0 \, \Delta f \, \omega
$$
Il résulte alors de la propriété 1 que
$$
f_r (0)
=
\frac{\ln (r)}{2\pi} \int_{\CR (0,r)} *df
-\frac{1}{2\pi}
\int_{\C\setminus D (0,r)} G_0 \, \Delta f\, \omega
\, .
$$
On conclut en faisant tendre $r$ vers 0 puisque l'intégrale $\int_{\CR (0,r)} *df$ est en $O(r)$ et que $G_0 \, \Delta f$ est intégrable sur $\C$ comme l'est, localement, $G_0$.
\cqfd

\niR
Il existe bien sûr des résultats similaires pour des cercles centrés en n'importe quel point $z_0$ de $\C$.
On peut les déduire de ceux obtenus ci-dessus, en utilisant par exemple l'action transitive et isométrique des translations du plan euclidien. On détaillera plus loin cette manière de procéder dans le cadre de la géométrie hyperbolique sur le disque unité.

\section{Un peu de géométrie dans le disque de Poincaré}

\subsection{La métrique hyperbolique sur le disque unité}

On considère sur le disque unité $\D=\{z\in \C\, |\ |z|<1\}$ du plan complexe
la structure différentiable induite par celle de $\C$
et l'orientation donnée par la forme volume $\omega=dx\wedge dy$.
La métrique hyperbolique $\ghyp$ sur $\D$ est définie par:
$$
\ghyp_z \left( u,v \right)
=
\frac1{\left(1-|z|^2\right)^2} \, g_{z} (u,v)\quad\quad\mbox{pour tous $u,v$ dans $T_z \D$ et tout $z$ dans $\D$}\, ,
$$
où $g$ désigne la métrique euclidienne sur $\C$ que l'on restreint ici à $\D$.
En particulier, pour tout $z$ dans $\D$, le couple de vecteurs
$$
\left(
\left(1-|z|^2\right){\partial_x}_{|z}
,
\left(1-|z|^2\right){\partial_y}_{|z}
\right)
$$
constitue une base orthonormée directe de $T_z\D$.

\v2
Fixons $z$ dans $\D$.
Comme on l'a fait pour la métrique euclidienne,
on associe à n'importe quel élément $u$ de $T_z \D$, un élément $u^{\starhypb}$
de $T_z^* \D$
défini par:
$$
u^{\starhypb} (v)=\ghyp_z(u,v)
\quad\quad
\mbox{pour tout $v$ dans $T_z\D$}\, .
$$
De part la définition de $\ghyp$ en fonction de $g$,
il vient:
$u^{\starhypb}=\left(1-|z|^2\right)^{-2} u^*$.
On en déduit que la forme volume canonique
$\omegahyp$ pour la métrique hyperbolique $\ghyp$ sur $\D$
est donnée par:
$$
\omegahyp_{|z}
=
\left(\left(1-|z|^2\right){\partial_x}_{|z}\right)^{\starhypb}\wedge\left(\left(1-|z|^2\right){\partial_x}_{|z}\right)^{\starhypb}
=
\left(1-|z|^2\right)^{-2}\,\omega_{|z}
\quad\quad
\mbox{pour tout $z$ dans $\D$}\,  .
$$
Ainsi,
l'opérateur étoile de Hodge pour la métrique euclidienne
et
l'opérateur étoile de Hodge pour la métrique hyperbolique
coïncident sur $T_z^* \D$ puisqu'on a l'équivalence:
$$
u^*\wedge *\beta=\beta (u)\, \omega_{|z}
\quad
\Longleftrightarrow
\quad
u^{\starhypb}\wedge *\beta=\beta (u)\, \omegahyp_{|z}
$$
pour tout $u$ dans $T_z \D$ et tout $\beta$ dans $T_z^* \D$.

\v2
Considérons à présent une fonction $f$ régulière sur un ouvert de $\D$.
Sa différentielle $df$ donne par dualité un champ de vecteurs
$\nablahyp f$, le gradient hyperbolique de $f$, caractérisé par:
$$
df_{|z} (v)=\ghyp_z (\nablahyp f_{|z},v) \quad \mbox{pour tout $v$ dans $T_z\D$ et tout $z$ dans l'ouvert considéré}\, .
$$
On vérifie que l'on a: $\nablahyp f_{|z}=\left(1-|z|^2\right)^{2} \nabla f_{|z}$.
De plus, comme l'opérateur étoile de Hodge ne diffère pas entre les métriques euclidienne et hyperbolique,
le laplacien hyperbolique $\Deltahyp f$ de $f$
vérifie:
$$
\Deltahyp f \, \omegahyp=d *\!df=\Delta f \, \omega\, ,
$$
et donc
$$
\Deltahyp f (z)=\left(1-|z|^2\right)^{2} \Delta f (z)
\ \ \mbox{pour tout $z$ dans l'ouvert considéré}\, .
$$
En particulier, une fonction harmonique sur un ouvert de $\D$ pour la métrique euclidienne l'est aussi pour la métrique hyperbolique, et réciproquement.

\subsection{Les isométries du disque de Poincaré}

Le disque $\D$ muni de la métrique hyperbolique est appelé disque de Poincaré.
Une isométrie du disque de Poincaré est un $C^\infty$-difféomorphisme $\Phi:\D\to \D$ qui préserve la métrique hyperbolique au sens où:
$$
\ghyp_{\Phi (z)} (T_z \Phi (u),T_z \Phi (v))=\ghyp_z (u,v)
$$
pour tous $u$ et $v$ dans $T_z \D$ et tout $z$ dans $\D$,
l'application $T_z \Phi:T_z \D \to T_{\Phi (z)} \D$ désignant l'application tangente à $\Phi$ en $z$.

\v2
La distance hyperbolique $\dhyp$ entre deux points $z$ et $z'$ de $\D$ est l'infimum des longueurs des chemins reliant $z$ à $z'$:
$$
\dhyp (z,z')
:=
\inf_{\gamma} \int_a^b \sqrt{\ghyp_{\gamma (t)} \left( \dot \gamma (t),\dot \gamma (t) \right)}\, dt \, ,
$$
où l'infimum est pris sur l'ensemble des chemins réguliers
$\gamma:[a,b]\to\D$ tels que $\gamma(a)=z$ et $\gamma (b)=z'$.
Une isométrie du disque de Poincaré préserve la distance hyperbolique $\dhyp$
puisqu'elle envoie l'ensemble des chemins réguliers joignant deux points de $\D$ sur l'ensemble des chemins réguliers joignant leurs images en préservant leur longueur comme on le voit grâce au théorème de dérivation des applications composées.

\v2
Une isométrie $\Phi$ du disque de Poincaré préserve également, au signe près, la forme volume canonique $\omegahyp$.
En effet, pour $z$ quelconque dans $\D$,
si  $(e_1,e_2)$ désigne une base orthonormée directe de $T_z \D$
alors $(T_z \Phi (e_1),T_z \Phi (e_2))$
est une base orthonormée de $T_{\Phi (z)} \D$ directe ou indirecte selon que $\Phi$ préserve l'orientation ou non,
et l'on a,
pour tous $u$ et $v$ dans $T_z\D$,
\begin{eqnarray*}
\left(\Phi^*\omegahyp\right)_{|z} (u\wedge v)
&=&
\vareps(\Phi)
\left( (T_z \Phi (e_1))^{*_\hyp} \wedge (T_z \Phi (e_2))^{*_\hyp} \right)_{|\Phi (z)}
\left(T_z \Phi (u)\wedge T_z \Phi (v)\right)
\\
&=&
\vareps(\Phi)
\left|
\begin{array}{cc}
(T_z \Phi (e_1))^{*_\hyp} (T_z \Phi (u)) & (T_z \Phi (e_1))^{*_\hyp} (T_z \Phi (v))\\
(T_z \Phi (e_2))^{*_\hyp} (T_z \Phi (u)) & (T_z \Phi (e_2))^{*_\hyp} (T_z \Phi (v))
\end{array}
\right|
\\
&=&
\vareps(\Phi)
\left|
\begin{array}{cc}
\ghyp_{\Phi (z)} (T_z \Phi (e_1), T_z \Phi (u)) & \ghyp_{\Phi (z)} ((T_z \Phi (e_1), T_z \Phi (v))\\
\ghyp_{\Phi (z)} (T_z \Phi (e_2), T_z \Phi (u)) & \ghyp_{\Phi (z)} ((T_z \Phi (e_2), T_z \Phi (v))
\end{array}
\right|
\\
&=&
\vareps(\Phi)
\left|
\begin{array}{cc}
\ghyp_{z} (e_1, u) & \ghyp_{z} (e_1,v)\\
\ghyp_{z} (e_2, u) & \ghyp_{z} (e_2,v)
\end{array}
\right|
=
\vareps(\Phi)\, \omegahyp_{|z} (u\wedge v)
\end{eqnarray*}
où $\vareps (\Phi)$ vaut $+1$ ou $-1$ selon que $\Phi$ préserve l'orientation ou non. Il en résulte qu'une isométrie $\Phi$ préserve l'aire hyperbolique puisque la formule de changement
de variables dans les intégrales de formes volumes donne, pour toute fonction $f$ positive et mesurable sur $\D$,
$$
\int_\D f\circ \Phi \ \omegahyp
=
\vareps (\Phi) \int_\D \Phi^* f \ \Phi^*\omegahyp
=
\vareps (\Phi) \int_\D \Phi^* (f \omegahyp)
=
\vareps (\Phi)^2 \int_\D f \omegahyp
=
\int_\D f \omegahyp
\, .
$$
Il en résulte également la relation de commutation: $*\, \Phi^*=\vareps (\Phi)\Phi^* *$, entre les opérateurs $\Phi^*$
et étoile de Hodge agissant sur les formes différentielles de degré 1 régulières sur un ouvert de $\C$. Comme de plus $\Phi^*$ et $d$ commutent,
on en déduit que,
pour toute fonction $f$ régulière sur un ouvert de $\D$, on a:
\begin{eqnarray*}
\Deltahyp (f\circ \Phi)\ \omegahyp
&=&
d*\!d(f\circ\Phi)
=
d*\!d (\Phi^* f)
=
d*\!(\Phi^* df)
\\
&=&
\vareps (\Phi)\
d\, \Phi^*(*df)
=
\vareps (\Phi)\
\Phi^* (d*\!d f)
=
\vareps (\Phi)\
\Phi^* (\Deltahyp f \, \omegahyp)
\\
&=&
\vareps (\Phi)\
\Phi^* (\Deltahyp f) \ \Phi^* \omegahyp
=
\vareps (\Phi)^2\
(\Deltahyp f)\circ \Phi \ \omegahyp
=
(\Deltahyp f)\circ \Phi \ \omegahyp\; ;
\end{eqnarray*}
ce qui montre que laplacien et isométrie commutent également.

\v2
Justifions pour terminer cette section que le groupe des isométries du disque de Poincaré agit transitivement sur $\D$.
\'Etant donné un point $z_0$ de $\D$
distinct de 0, il existe une unique inversion $\Phi_{z_0}$ du plan complexe euclidien laissant globalement invariant $\D$ et envoyant $z_0$ sur 0. Elle est donnée par l'anti-homographie (involutive):
$$
\Phi_{z_0} (z)=\frac{z_0}{\o{z_0}} \frac{\o{z}-\o{z_0}}{z_0 \o{z}-1}\, .
$$
La restriction de $\Phi_{z_0}$ à $\D$, que l'on désigne encore par $\Phi_{z_0}$, est une isométrie du disque de Poincaré car c'est le cas de la conjugaison $z\mapsto \o{z}$ comme on le vérifie facilement,
ainsi que de l'homographie $z\mapsto \o{\Phi}_{z_0} (z):=\o{\Phi_{z_0} (z)}$
puisque l'on a\footnote{utiliser de plus que si l'on identifie $T_z \D$ à $\C$ via la base
$({\partial_x}_{|z},{\partial_y}_{|z})$
alors on a: $T_z \o{\Phi}_{z_0} (u)=\o{\Phi}_{z_0}' (z) u$, et $\ghyp_{z} (u,v)=(1-|z|^2)^{-2}\frac{\o{u}v+u\o{v}}{2}$}:
$$
\frac{\left|\o{\Phi}_{z_0}' (z)\right|^2}{\left(1-\left|\o{\Phi}_{z_0} (z)\right|^2\right)^2}=\frac{1}{\left(1-|z|^2\right)^2}\, .
$$
On en déduit la transitivité de l'action sur $\D$ du groupe des isométries du disque de Poincaré\footnote{On peut montrer plus précisément que les isométries du disque de Poincaré sont exactement les homographies
et anti-homographies de $\C$ qui préservent $\D$ et que leur groupe est engendré par les inversions de $\C$ qui préserve $\D$ (voir \cite{Beardon} par exemple)}.

\section{Moyenner sur les cercles les fonctions régulières dans le disque de Poincaré}

On suppose le disque $\D$ muni de la métrique hyperbolique $\ghyp$,
de la forme volume canonique $\omegahyp$
et de la distance $\dhyp$.

\v2
Puisque la métrique hyperbolique sur $\D$ se déduit de la métrique euclidienne via une fonction à symétrie radiale,
pour tout point $z$ de $\D$,
le segment euclidien $[0,z]$ est aussi le plus court chemin de 0 à $z$ pour la métrique hyperbolique
et l'on a:
$$
\dhyp (0,z)
=
\int_0^{|z|} \frac{dr}{1-r^2}
=
\frac12 \ln \left(\frac{1+|z|}{1-|z|} \right)
\, ;
\quad\quad
\mbox{ou encore}\, ,
\quad
|z| = \tanh \left( \dhyp (0,z) \right)\, .
$$
Ainsi, pour $r=\tanh ( \rho )$ avec $\rho>0$,
le cercle $\CR^{\hyp} (0,\rho)$ centré en 0 et de rayon $\rho$
du disque de Poincaré
coïncide en tant qu'ensemble
avec le cercle $\CR (0,r)$ centré en 0 et de rayon $r$ du disque euclidien.
De plus,
toujours du fait que les deux métriques considérées sont conformément équivalentes via une fonction à symétrie radiale, pour toute fonction $f$ mesurable sur $\CR (0,r)$,
la moyenne hyperbolique
$f_{\rho}^{\hyp} (0)$ de $f$
sur $\CR^{\hyp} (0,\rho)$\footnote{c'est-à-dire la moyenne de $f$ sur $\CR^{\hyp} (0,\rho)$ lorsque ce dernier est muni de la métrique induite par la métrique hyperbolique sur $\D$}
coïncide
avec 
la moyenne euclidienne $f_r (0)$ de $f$ sur $\CR (0,r)$:
$$
f_{\rho}^{\hyp} (0) = f_r (0)\quad\quad \mbox{pour $r=\tanh ( \rho )$ avec $\rho>0$}\, .
$$

\v2
Considérons à présent un point $z_0$ de $\D$
et
une isométrie $\Phi_{z_0}$ du disque de Poincaré
qui envoie $z_0$ sur 0\footnote{par exemple l'inversion $\Phi_{z_0}$ de la section précédente}.
\'Etant donnée une fonction $f$ mesurable sur le cercle $\CR^{\hyp} (z_0,\rho)$ de centre $z_0$ et de rayon $\rho>0$ pour la métrique hyperbolique,
la moyenne hyperbolique $f_{\rho}^{\hyp} (z_0)$ de $f$ sur $\CR^{\hyp} (z_0,\rho)$
coïncide avec
la moyenne hyperbolique $\left( f\circ \Phi_{z_0}^{-1}\right)_{\rho}^{\hyp} (0)$ de $f\circ \Phi_{z_0}^{-1}$ sur $\CR^{\hyp} (0,\rho)$.
On a donc:
$$
f_{\rho}^{\hyp} (z_0)
=
\left( f\circ \Phi_{z_0}^{-1}\right)_{\rho}^{\hyp} (0)=
\left( f\circ \Phi_{z_0}^{-1} \right)_r (0) \quad\quad \mbox{pour $r=\tanh ( \rho )$}\, .
$$

\v2
Le résultat suivant se déduit à présent de la propriété 2 ci-dessus.

\v1
\n{\bf Propriété 4 --}
Si $f$ est une fonction régulière et harmonique
sur le disque ouvert $D^{\hyp} (z_0,\rho)$
de centre $z_0$ et de rayon $\rho>0$ pour la métrique hyperbolique
et si $f$ est continue sur l'adhérence de $D^{\hyp} (z_0,\rho)$
alors on a:
$$
f (z_0)=f^{\hyp}_\rho (z_0)\, .
$$

\niP
En effet,
avec les notations précédemment introduites, il vient:
$$
f (z_0)=f\circ \Phi_{z_0}^{-1} (0) 
=
\left( f\circ \Phi_{z_0}^{-1} \right)_r (0)
=
\left( f \circ \Phi_{z_0}^{-1} \right)^{\hyp}_\rho (0)
=
f^{\hyp}_\rho (z_0)
$$
car $f\circ \Phi_{z_0}^{-1}$ est harmonique (pour les deux métriques) sur $D(0,r)$ puisque
$\Phi_{z_0}^{-1}$ et $\Deltahyp$ commutent.

\v2
La propriété suivante
se déduit de la propriété 3 ci-dessus.

\v1
\n{\bf Propriété 5 --}
Si $f$ est une fonction définie et régulière sur $\D$, à support compact,
alors, pour tout $z_0$ dans $\D$, on a:
$$
f (z_0)
=
-\frac{1}{2\pi}
\int_{\D} G_\D (z_0,\cdot)  \, \Deltahyp f \, \omegahyp
=
-\frac{1}{2\pi}
\int_{\D} G_\D (z_0,\cdot)  \, \Delta f \, \omega
$$
où l'on a introduit la fonction de Green
$G_\D$ de $\D$
définie par
$$
G_\D (z_0,z)=
-\ln \left( \tanh \left( \dhyp (z_0,z) \right) \right)
\quad \mbox{pour tout $z$ dans $\D\setminus\{z_0\}$}\, .
$$

\niP
En effet,
on a, d'une part,
\begin{eqnarray*}
f(z_0)
&=&
f\circ \Phi_{z_0}^{-1} (0)
=
-\frac{1}{2\pi}
\int_{\D} G_0\ \Delta \left(f\circ\Phi_{z_0}^{-1}\right) \ \omega
\\
&=&
-\frac{1}{2\pi}
\int_{\D} G_0\ \Deltahyp \left(f\circ\Phi_{z_0}^{-1}\right) \ \omegahyp
=
-\frac{1}{2\pi}
\int_{\D} G_0\ \left(\Deltahyp f\right)\circ\Phi_{z_0}^{-1} \ \omegahyp
\\
&=&
-\frac{1}{2\pi}
\int_{\D} G_0\circ\Phi_{z_0}\ \Deltahyp f \ \omegahyp
=
-\frac{1}{2\pi}
\int_{\D} G_0\circ\Phi_{z_0}\ \Delta f \ \omega_{0}
\end{eqnarray*}
et d'autre part, pour tout $z$ dans $\D$,
\begin{eqnarray*}
G_0\circ\Phi_{z_0} (z)
&=&
-
\ln \left( \left|\Phi_{z_0} (z) \right| \right)
=
-
\ln \left( \tanh \left( \dhyp (0,\Phi_{z_0} (z)) \right) \right)
\\
&=&
-
\ln \left( \tanh \left( \dhyp (\Phi_{z_0}^{-1} (0),z) \right) \right)
=
-
\ln \left( \tanh \left( \dhyp (z_0,z) \right) \right)\, .
\end{eqnarray*}
\cqfd

\v5
Pour terminer cette section et introduire la suivante, remarquons qu'avec les notations précédentes, pour toute fonction $f$ définie et régulière sur $\D$, à support compact, on a:
\begin{eqnarray*}
f_\rho^{\hyp} (z_0)
&=&
\left(f \circ \Phi_{z_0}^{-1}\right)_\rho^{\hyp} (0)
=
\left(f \circ \Phi_{z_0}^{-1}\right)_r (0)
\\
&=&
\frac{1}{2\pi}
\int_{\D\setminus D (0,r)}
g\left(\nabla G_0, \nabla \left(f\circ \Phi_{z_0}^{-1}\right)\right) \, \omega
\end{eqnarray*}
avec $\rho>0$ et $r=\tanh (\rho)$.
Ainsi, en introduisant la fonction
$G_0^r$ définie
par:
$$
G_0^r(z)
=
-
\ln \left(\max(r,|z|) \right)\quad \mbox{pour tout $z$ dans $\D$}\, ,
$$
il vient:
$$
f_\rho^{\hyp} (z_0)
=
\frac{1}{2\pi}
\int_{\D}
g\left(\nabla G_0^r, \nabla \left(f\circ \Phi_{z_0}^{-1}\right)\right) \, \omega
\, ,
$$
où $\nabla G_0^r$ désigne le gradient (euclidien) au sens des distributions de la fonction $G_0^r$:
$$
\nabla G_0^r
=
-
\left(\frac{x}{x^2+y^2}\, {\partial_x}
+
\frac{y}{x^2+y^2}\, {\partial_y} \right)\,  \1_{\{r<|z|<1\}}\, .
$$

\section{L'espace de Sobolev $H_0^1 (\D)$}

On munit
l'espace vectoriel $C_c^\infty (\D)$ des fonctions définies et régulières sur $\D$,
à support compact,
d'une structure d'espace préhilbertien grâce au produit scalaire défini par:
$$
\la \varphi,\psi \ra
=
\frac{1}{2\pi}
\int_\D g(\nabla \varphi ,\nabla \psi)\, \omega
\quad\quad\mbox{pour tous $\varphi$ et $\psi$ dans $C^\infty (\D)$}\, .
$$
(Un élément de ``petite'' norme dans cet espace est donc une fonction présentant de ``petites'' variations.)

Notons dès à présent que cet espace préhilbertien ne dépend pas de la métrique euclidienne ou hyperbolique sur $\D$ puisque l'on a:
$$
\la \varphi,\psi \ra
=
\frac{1}{2\pi}
\int_\D d\varphi \wedge *d\psi
=
\frac{1}{2\pi}
\int_\D \ghyp (\nablahyp \varphi ,\nablahyp \psi)\, \omegahyp\, .
$$

Pour obtenir une représentation pratique du complété de $(C_c^\infty (\D),\la \cdot,\cdot \ra)$,
on peut procéder comme suit.
Rappelons\footnote{On renvoie à \cite{Brezis} pour une présentation systématique de ces questions.} que l'espace $H^1 (\D)$ désigne l'espace de Hilbert (séparable)
des fonctions de $L^2 (\D)$, dont les dérivées partielles d'ordre 1 au sens des distributions
sont également
dans $L^2 (\D)$, muni du produit scalaire défini par:
$$
\la u,v \ra_{H^1 (\D)}
=
\int_\D u v\, \omega
+
\int_\D g(\nabla u ,\nabla v)\, \omega
\quad\quad\mbox{pour tous $u$ et $v$ dans $H^1 (\D)$}\, .
$$
L'espace vectoriel $H_0^1 (\D)$ désigne quant à lui l'adhérence
de $C_c^\infty (\D)$ dans $H^1 (\D)$.
L'inégalité de Poincaré permet de justifier que l'on définit
un nouveau produit scalaire sur $H_0^1 (\D)$ en posant:
$$
\la u,v \ra
=
\frac{1}{2\pi}
\int_\D g(\nabla u ,\nabla v)\, \omega
\quad\quad\mbox{pour tous $u$ et $v$ dans $H_0^1 (\D)$}\, ,
$$
et que les deux normes hilbertiennes ainsi obtenues sur $H_0^1 (\D)$ sont équivalentes.
Avec ce nouveau produit scalaire, $H_0^1 (\D)$ est donc encore un espace de Hilbert
qui donne bien une représentation du complété de $(C_c^\infty (\D),\la \cdot,\cdot \ra)$.
C'est de cette représentation dont il s'agit quand on évoque $H_0^1 (\D)$
dans la suite de cette note.

\v3
La propriété suivante permet d'identifier une classe importante d'éléments de $H_0^1 (\D)$.
On en donne une preuve succincte
car ses ingrédients seront utilisés dans la prochaine section.

\v1
\n{\bf Propriété 6 --}
Soit $f$ une fonction continue sur $\D$ dont le prolongement $\tilde{f}$ par 0 sur $\C\setminus \D$
est continu sur $\C$.
Si, de plus, $f$ est élément de $H^1 (\D)$ alors $f$ appartient à $H_0^1 (\D)$.

\niP
Il s'agit de montrer que la fonction $f$ peut être approchée d'aussi près que l'on veut dans $H^1 (\D)$
par des éléments de $C_c^\infty (\D)$.

On le vérifie lorsque $f$ est à support compact (inclus dans $\D$) en considérant une suite régularisante 
$(\rho_n)_{n\geq 1}$ de fonctions définies sur $\C$ telle que, pour tout entier $n\geq 1$,
la restriction de $\rho_n*\tilde{f}$ à $\D$ appartienne à $C_c^\infty (\D)$,
et telle que l'on ait:
$$
\lim_{n\to +\infty} \sup_{z\in \D} \left|\rho_n*\tilde{f} (z) - f(z)\right|=0
$$
et
\begin{eqnarray*}
& &
\int_\D g\left(\nabla (\rho_n*\tilde{f}-f),\nabla (\rho_n*\tilde{f} - f)\right)\, \omega
\\
&=&
\int_\D \left|\rho_n * (\partial_x \tilde{f}) - \partial_x f\right|^2 \, \omega
+
\int_\D \left|\rho_n * (\partial_y \tilde{f}) - \partial_y f\right|^2 \, \omega
\quad
\underset{n\to +\infty}{\longrightarrow}
\quad 0\, .
\end{eqnarray*}

Pour conclure dans le cas général,
il suffit donc de montrer que l'on peut approcher dans $H^1 (\D)$ la fonction $f$ par une suite $(f_n)_{n\geq 1}$ de fonctions
continues sur $\D$ à support compact (inclus dans $\D$).
En fait, la suite de fonctions définie par
$$
f_n (z)=\frac{1}{n} h(nf(z))\, , \quad \mbox{ pour tout $z$ dans $\D$ et tout $n\geq 1$}\, ,
$$
convient dès que $h$ désigne une fonction de classe $C^1$ sur $\R$ qui vérifie, pour tout $t$ dans $\R$,
$$
|h(t)|\leq |t|
\qeq
h(t)=
\left\{
\begin{array}{ll}
0 & \mbox{si $|t|\leq 1$}\\
t & \mbox{si $|t|\geq 2$}
\end{array}
\right.
$$
En effet,
pour tout $n\geq 1$,
le support de $f_n$ est alors contenu dans $\{z\in \D\, |\ |f(z)|\geq 1/n\}$ et est donc, par hypothèse sur $f$, un compact inclus dans $\D$. De plus, pour tout $n\geq 1$, $f_n$ est continue sur $\D$ comme composée de fonctions continues et l'on a:
$$
\sup_{z\in \D} \left|f_n (z) - f(z)\right|\leq \frac{4}{n}\, .
$$
Enfin,
comme $h(0)=0$ et comme la dérivée $h'$ de $h$ est bornée sur $\R$,
la fonction $f_n$ est élément de $H^1 (\D)$ et il vient, pour tout $n\geq 1$:
$$
\nabla f_n =h'(nf)\, \nabla f\quad\quad\mbox{au sens des distributions}\, .
$$
On en déduit la convergence:
$$
\nabla f_n
\ 
\underset{n\to +\infty}{\longrightarrow}
\ 
\nabla f
\quad\mbox{presque partout}
$$
(utiliser que $\nabla f=0$ presque partout sur $\{f=0\}$\footnote{voir \cite{Gilbarg-Trudinger} p. 145})
et le théorème de convergence dominée permet de conclure.
\cqfd

\niE
La propriété 6 permet en particulier de justifier l'appartenance à $H_0^1 (\D)$
de la fonction $G_0^r$ introduite en fin de section précédente.

\v5
Terminons cette section en remarquant que,
pour toutes fonctions $\varphi$ et $\psi$ dans $C_c^\infty (\D)$,
on a:
$$
\la \varphi,\psi \ra
=
- \frac{1}{2\pi} \int_\D \varphi\, \Delta \psi\, \omega
=
- \frac{1}{2\pi} \int_\D \varphi\, \Deltahyp \psi\, \omegahyp\, ;
$$
et qu'ainsi,
pour toute une isométrie $\Phi$ du disque de Poincaré, il vient:
\begin{eqnarray*}
\la \varphi\circ \Phi,\psi\circ \Phi\ra
&=&
- \frac{1}{2\pi} \int_\D \varphi\circ \Phi\, \Deltahyp (\psi\circ \Phi)\, \omegahyp
=
- \frac{1}{2\pi} \int_\D \varphi\circ \Phi\, \left(\Deltahyp \psi\right)\circ \Phi\ \omegahyp
\\
&=&
- \frac{1}{2\pi} \int_\D \varphi\, \Deltahyp \psi \, \omegahyp
=
\la \varphi,\psi \ra
\, ,
\end{eqnarray*}
où l'on a utilisé que
les fonctions $\varphi\circ \Phi$ et $\psi\circ \Phi$ sont dans $C_c^\infty (\D)$,
que la composition par $\Phi$ et le laplacien hyperbolique commutent
et que $\Phi$ préserve l'aire hyperbolique.
Ainsi, l'application $\Phi^*:\varphi\mapsto \varphi\circ \Phi$ est une isométrie linéaire de l'espace préhilbertien $C_c^\infty (\D)$. Elle se prolonge par densité en une isométrie linéaire $\Phi^*$ de $H_0^1 (\D)$ et l'on a encore l'égalité:
$\Phi^* (u)=u\circ \Phi$ pour tout $u$ dans $H_0^1 (\D)$
(la convergence d'une suite dans $L^2$ assurant sa convergence presque partout le long d'une sous-suite).

\section{Moyenner sur les cercles les éléments de $H_0^1 (\D)$ dans le disque de Poincaré}

Fixons un point $z_0$ de $\D$ et $\Phi_{z_0}$ une isométrie du disque de Poincaré qui envoie $z_0$ sur 0.
Considérons également un  nombre réel $\rho$ strictement positif et posons
$r=\tanh \left( \rho \right)$.

On a remarqué en fin de section 5 que
la moyenne hyperbolique sur le cercle $\CR^\hyp (z_0,\rho)$ d'un élément $\varphi$
de $C_c^\infty (\D)$ vérifie l'égalité:
$$
\varphi_\rho^{\hyp} (z_0)
=
\frac{1}{2\pi}
\int_{\D}
g\left(\nabla G_0^r, \nabla \left(\varphi\circ \Phi_{z_0}^{-1}\right)\right) \, \omega\, ,
$$
avec
$G_0^r(z)=- \ln \left(\max(r,|z|) \right)$.
Comme $G_0^r$ est élément de $H_0^1 (\D)$, on peut aussi écrire:
$$
\varphi_\rho^{\hyp} (z_0)
=
\la G_0^r , \varphi\circ \Phi_{z_0}^{-1} \ra
\, ;
$$
et donc également, en utilisant l'invariance du produit scalaire sous l'action des isométries du disque de Poincaré,
$$
\varphi_\rho^{\hyp} (z_0)
=
\la G_0^r \circ \Phi_{z_0},\varphi \ra
\, .
$$

\v5
La propriété suivante étend cette dernière égalité aux fonctions continues de $H_0^1 (\D)$
qui ont été considérées dans la propriété 6.

\v1
\n{\bf Propriété 7 --}
Soit $f$ un élément de $H^1 (\D)$.
Soient $z_0$ un point de $\D$ et $\rho$ un nombre réel strictement positif.

Si $f$ est une fonction continue sur $\D$ et si son prolongement par 0 sur $\C\setminus \D$ est continu sur $\C$
alors
on a:
$$
f_\rho^{\hyp} (z_0) = \la G_0^r\circ \Phi_{z_0},f \ra\quad\quad\mbox{où $r=\tanh \left( \rho \right)$}\, .
$$

\niP
C'est une conséquence du fait que, comme il a été vu dans la preuve de la propriété 6,
on peut approcher $f$ dans $H_0^1 (\D)$ par une suite d'éléments de $C_c^\infty (\D)$ qui, en outre, converge uniformément sur $\D$ vers $f$
\cqfd

\niR
La théorie des traces permet de généraliser la propriété précédente à tout
élément de $H_0^1 (\D)$.

\section{Le champ libre gaussien sur le disque unité}

Commençons par une petite digression en dimension finie.
Un vecteur gaussien standard de $\R^d$ est un vecteur aléatoire $W=(W_1,\ldots,W_d)$ dont les coordonnées sont indépendantes et de même loi gaussienne centrée de variance unité. D'un point de vue plus géométrique,
un vecteur gaussien standard de $\R^d$ est, à un scalaire près,
un vecteur gaussien qui admet une densité
dont les (hyper)surfaces de niveaux sont les sphères euclidiennes de $\R^d$;
autrement dit, un vecteur gaussien qui ne privilégie aucune direction particulière.
Plus généralement,
une variable aléatoire gaussienne standard à valeurs dans un espace euclidien $(E,\la \cdot,\cdot \ra)$
est une variable aléatoire $W$
à valeurs dans $E$
telle que, pour tout vecteur $u$ de $E$, la variable aléatoire réelle $\la u,W \ra$ est gaussienne, centrée et de variance $\norm{u}{}^2$.

\v2
Alors que l'existence d'une variable gaussienne standard à valeurs dans n'importe quel espace euclidien
est assurée, il n'existe pas de variable aléatoire $W$ à valeurs dans un espace de Hilbert réel
$(H,\la \cdot,\cdot \ra)$ de dimension infinie pour laquelle
la variable aléatoire réelle $\la u,W \ra$ serait gaussienne, centrée et de variance $\norm{u}{}^2$
pour n'importe quel vecteur $u$ dans $H$
(raisonner par l'absurde,
considérer une famille orthonormée $(e_i)_{i\in \N}$ de $H$
et remarquer que la suite $(\la e_i,W \ra)_{i\in \N}$ de variables aléatoires réelles
de même loi gaussienne standard
devrait converger presque sûrement vers 0, ce qui est exclus).
\`A défaut,
on peut
quand même
construire une isométrie linéaire $W$
de $H$ dans un espace $L^2 (\Omega,\FR,\P)$ 
telle que l'image $W(u)$ de n'importe quel vecteur $u$ de $H$
soit une variable aléatoire réelle, gaussienne et centrée (de variance nécessairement égale à $\norm{u}{}^2$).
En effet, il suffit de considérer une base hilbertienne $(e_i)_{i\in I}$ de $H$ et
un espace $L^2 (\Omega,\FR,\P)$ sur lequel il existe une famille $(W_i)_{i\in I}$
de variables aléatoires réelles indépendantes de même loi gaussienne standard\footnote{un tel espace existe pour n'importe quel ensemble $I$ par le théorème de prolongement de Kolmogorov (voir par exemple \cite{Neveu} p. 79)};
l'application $W$ qui, à tout vecteur $u$ de $H$, associe la série $\sum_{i\in I} \la e_i,u\ra W_i$
convient\footnote{la sommabilité de la série $\sum_{i\in I} \la e_i,u\ra W_i$ dans $L^2 (\Omega,\FR,\P)$
est une conséquence de la sommabilité de la série $\sum_{i\in I} \la e_i,u\ra$ dans $\ell^2 (I)$,
du critère de Cauchy et de l'orthonormalité des $W_i$}.
La variable aléatoire $W(u)$ est souvent notée de manière
abusive (mais allusive) $\la u,W \ra$.

\v2
Une telle isométrie $W$, lorsqu'elle est définie sur l'espace de Hilbert $H_0^1 (\D)$,
est
appelée
champ libre gaussien sur le disque unité $\D$.

\niR
\'Etant donné un champ libre gaussien $W$ sur $\D$,
on a,
en utilisant la propriété 5,
pour toutes fonctions $\varphi_1$ et $\varphi_2$ dans $C_c^\infty (\D)$:
\begin{eqnarray*}
\cov(\la\varphi_1,W\ra,\la\varphi_2,W\ra)
&=&
\la \varphi_{1},\varphi_{2} \ra
=
\frac{1}{2\pi}
\int_\D (-\Delta \varphi_{1})\, \varphi_{2}\, \omega
\\
&=&
\frac{1}{4\pi^2}
\int_{\D\times\D} (-\Delta \varphi_{1}) (z_1)\, (-\Delta \varphi_{2}) (z_2)\, G_\D (z_1,z_2) \, \omega\otimes\omega_{|(z_1,z_2)}\, .
\end{eqnarray*}
On peut rapprocher cette expression de la formule de covariance d'un pont brownien standard $(Z_t)_{t\in [0,1]}$:
$$
\cov (Z_s, Z_t)=\min(s,t)(1-\max(s,t))=:G_{]0,1[} (s,t)\, ,
$$
après avoir remarqué que, pour toute fonction $\varphi$ définie et régulière sur $]0,1[$,
à support compact,
on a:
$$
\varphi (t)=\int_0^1 G_{]0,1[} (s,t)\, \left(-\varphi''\right) (s)\, ds\, ,
\quad
\mbox{pour tout $t$ dans $]0,1[$}\, .
$$

\section{Moyenner sur les cercles le champ libre gaussien dans le disque de Poincaré}

Considérons un champ libre gaussien $W$ sur $\D$.

\'Etant donnés un point $z$ de $\D$ et un nombre réel $\rho$ strictement positif,
les considérations de la section 7 invitent à appeler moyenne hyperbolique de $W$ sur le cercle $\CR^{\hyp} (z,\rho)$
la variable aléatoire:
$$
W_\rho^\hyp (z) := \la G_0^r\circ \Phi_{z},W \ra\quad\quad
\mbox{où $r=\tanh(\rho)$}\, .
$$

Par définition de $W$,
le processus $(W_\rho^\hyp (z)\, |\ (z,\rho)\in \D\times ]0,+\infty[)$ est un processus gaussien centré dont la fonction de covariance est donnée par:
$$
\left( (z_1,\rho_1),(z_2,\rho_2) \right)
\longmapsto
\cov\left( W^\hyp_{\rho_1} (z_1) , W^\hyp_{\rho_2} (z_2) \right)
=
\la G_0^{r_1} \circ \Phi_{z_1} , G_0^{r_2} \circ \Phi_{z_2}\ra\, ,
$$
avec $r_1=\tanh (\rho_1)$ et $r_2=\tanh (\rho_2)$.
Ainsi, par la propriété 7,
la covariance
entre $W^\hyp_{\rho_1} (z_1)$ et $W^\hyp_{\rho_2} (z_2)$
est égale
à la moyenne hyperbolique de la fonction $G_0^{r_2} \circ \Phi_{z_2}$ sur $\CR^{\hyp} (z_1,\rho_1)$,
mais aussi à la moyenne hyperbolique de $G_0^{r_1} \circ \Phi_{z_1}$ sur $\CR^{\hyp} (z_2,\rho_2)$.
La propriété et les deux théorèmes qui suivent en sont des conséquences.

\v5
\n{\bf Propriété 8 --}
Soient $z_1$ et $z_2$ dans $\D$ et $\rho_1,\rho_2>0$.

Lorsque $\dhyp (z_1,z_2) \leq |\rho_2-\rho_1|$, on a:
$$
\cov\left( W^\hyp_{\rho_1} (z_1) , W^\hyp_{\rho_2} (z_2) \right)
=
-\ln \left(\tanh (\max(\rho_1,\rho_2))\right)\, .
$$

En particulier, lorsque $z_1=z_2=:z$ et $\rho_1=\rho_2=:\rho$, il vient:
$$
\var (W_{\rho}^\hyp (z))
=
-\ln \left(\tanh (\rho)\right)\, .
$$

\niP
Supposons par exemple $\rho_2$ supérieur à $\rho_1$. Alors, par hypothèse,
le disque intérieur au cercle $\CR^{\hyp} (z_2,\rho_2)$ contient le cercle $\CR^{\hyp} (z_1,\rho_1)$.
La fonction $G_0^{r_2} \circ \Phi_{z_2}$ est donc constante égale à $-\ln (\tanh (\rho_2))$ sur $\CR^{\hyp} (z_1,\rho_1)$; ce qui permet de conclure.\cqfd

\v5
\n{\bf Théorème 1 --}
Soient $z_1$ et $z_2$ dans $\D$ et $\rho_1,\rho_2>0$.

Lorsque $\dhyp (z_1,z_2)\geq\rho_1+\rho_2$, on a:
$$
\cov\left( W^\hyp_{\rho_1} (z_1) , W^\hyp_{\rho_2} (z_2) \right)
=
-\ln \left(\tanh \left( \dhyp (z_1,z_2) \right)\right)
= G_\D (z_1,z_2)
\, .
$$

\niP
Par hypothèse, les disques ouverts $D^{\hyp} (z_1,\rho_1)$ et $D^{\hyp} (z_2,\rho_2)$
intérieurs aux cercles $\CR^{\hyp} (z_1,\rho_1)$ et $\CR^{\hyp} (z_2,\rho_2)$ respectivement
sont disjoints. Ainsi,
la fonction $G_0^{r_1} \circ \Phi_{z_1}$ est harmonique sur le disque $D^{\hyp} (z_2,\rho_2)$ et continue sur son adhérence. On conclut alors en remarquant que, d'après la propriété 4, la moyenne hyperbolique
de la fonction $G_0^{r_1} \circ \Phi_{z_1}$ sur le cercle $\CR^{\hyp} (z_2,\rho_2)$ est égale à la valeur de $G_0^{r_1} \circ \Phi_{z_1}$ en le centre $z_2$ du cercle.
\cqfd

\niR
Ainsi, lorsque l'on a: $\dhyp (z_1,z_2)\geq \rho_1+\rho_2$ (c'est-à-dire pour $\rho_1$ et $\rho_2$ assez petits, $z_1$ et $z_2$ étant donnés), il résulte du théorème précédent que
la covariance entre $W^\hyp_{\rho_1} (z_1)$ et $W^\hyp_{\rho_2} (z_2)$ ne dépend que de la distance hyperbolique entre $z_1$ et $z_2$.

\v5
\n{\bf Théorème 2 --}
Pour tous $z_1$ et $z_2$ dans $\D$ et tous réels strictement positifs $\rho_1$ et $\rho_2$,
on a:
\begin{eqnarray*}
& &\E \left( \left( W^\hyp_{\rho_1} (z_1) - W^\hyp_{\rho_2} (z_2) \right)^2\right)\\
&\leq&
|\ln (\tanh (\rho_1)) - \ln(\tanh (\rho_2))|
+
2 \left( \frac{1}{\sinh(2\rho_1)} + \frac{1}{\sinh(2\rho_2)}\right)\, \dhyp (z_1,z_2)\, .
\end{eqnarray*}

\niP
Par définition du champ libre gaussien sur $\D$
et de ses moyennes hyperboliques sur les cercles, on a:
$$
\E \left( \left( W^\hyp_{\rho_1} (z_1) - W^\hyp_{\rho_2} (z_2) \right)^2\right)
=
\E \left(\la G_0^{r_1} \circ \Phi_{z_1} - G_0^{r_2} \circ \Phi_{z_2}, W\ra^2 \right) 
=
\norm{G_0^{r_1} \circ \Phi_{z_1} - G_0^{r_2} \circ \Phi_{z_2}}{}^2\, ,
$$
avec $r_1=\tanh (\rho_1)$ et $r_2=\tanh (\rho_2)$.
On en déduit, par invariance du produit scalaire sur $H_0^1 (\D)$
relativement à l'action des isométries du disque de Poincaré
sur $H_0^1 (\D)$, que
\begin{eqnarray*}
& &
\E \left( \left( W^\hyp_{\rho_1} (z_1) - W^\hyp_{\rho_2} (z_2) \right)^2\right)
\\
&=&
\norm{G_0^{r_1}}{}^2
+ \norm{G_0^{r_2}}{}^2
- 2 \la G_0^{r_1} , G_0^{r_2}\ra
\\
& & +
\la G_0^{r_1} \circ \Phi_{z_1} , G_0^{r_2} \circ \Phi_{z_1}\ra - \la G_0^{r_1} \circ \Phi_{z_1} , G_0^{r_2} \circ \Phi_{z_2}\ra
+
\la G_0^{r_1} \circ \Phi_{z_2} , G_0^{r_2} \circ \Phi_{z_2}\ra - \la G_0^{r_1} \circ \Phi_{z_1} , G_0^{r_2} \circ \Phi_{z_2}\ra
\\
&=&
\ln (\max(r_1,r_2)) - \ln (\min(r_1,r_2))
\\
& & +
\la G_0^{r_1} \circ \Phi_{z_1} ,G_0^{r_2} \circ \Phi_{z_1} - G_0^{r_2} \circ \Phi_{z_2}\ra
+
\la G_0^{r_1} \circ \Phi_{z_2}  - G_0^{r_1} \circ \Phi_{z_1} , G_0^{r_2} \circ \Phi_{z_2}\ra
\\
&\leq&
|\ln (r_1) - \ln(r_2)|
+
\left| \la G_0^{r_1} \circ \Phi_{z_1} , G_0^{r_2} \circ \Phi_{z_1} - G_0^{r_2} \circ \Phi_{z_2}\ra \right|
+
\left| \la G_0^{r_1} \circ \Phi_{z_2}  - G_0^{r_1} \circ \Phi_{z_1} , G_0^{r_2} \circ \Phi_{z_2}\ra \right|
\, .
\end{eqnarray*}

Or, la moyenne hyperbolique
$$
\la G_0^{r_1} \circ \Phi_{z_1} , G_0^{r_2} \circ \Phi_{z_1} - G_0^{r_2} \circ \Phi_{z_2}\ra
$$
de la fonction
$G_0^{r_2} \circ \Phi_{z_1} - G_0^{r_2} \circ \Phi_{z_2}$
sur le cercle $\CR^{\hyp} (z_1,\rho_1)$
est en valeur absolue majorée par la quantité:
$$
\sup_{z\in\D} \left|G_0^{r_2} \circ \Phi_{z_1} (z) - G_0^{r_2} \circ \Phi_{z_2} (z)\right|\, .
$$
Et comme, pour tout $z$ dans $\D$ et pour $i=1,2$, il vient:
$$
G_0^{r_2} \circ \Phi_{z_i} (z)=
-\ln \left(\max \left( r_2,\left|\Phi_{z_i} (z)\right| \right)\right)
=
-\ln \left(\tanh \left(\max \left( \rho_2,\dhyp (z,z_i) \right)\right)\right)\, ,
$$
il résulte de l'inégalité des accroissements finis que l'on a:
$$
\left|G_0^{r_2} \circ \Phi_{z_1} (z) - G_0^{r_2} \circ \Phi_{z_2} (z)\right|
\leq
\left( \frac{1}{\tanh (\rho_2)} - \tanh (\rho_2)\right) \left| \dhyp (z,z_1)-\dhyp (z,z_2)\right|
\, ;
$$
d'où l'on déduit l'inégalité:
$$
\left|
\la G_0^{r_1} \circ \Phi_{z_1} , G_0^{r_2} \circ \Phi_{z_1} - G_0^{r_2} \circ \Phi_{z_2} \ra
\right|
\leq \frac{2}{\sinh(2\rho_2)}\, \dhyp (z_1,z_2)
\, .
$$

On conclut en procédant de même avec la moyenne hyperbolique
$$
\la G_0^{r_1} \circ \Phi_{z_2}  - G_0^{r_1} \circ \Phi_{z_1} , G_0^{r_2} \circ \Phi_{z_2}\ra
$$
de $G_0^{r_1} \circ \Phi_{z_2}  - G_0^{r_1} \circ \Phi_{z_1}$ sur $\CR^{\hyp} (z_2,\rho_2)$.
\cqfd

\niR
La méthode de Kolmogorov-Chentsov
\footnote{voir
\cite{Revuz-Yor} et \cite{Khoshnevisan}
par exemple}
permet de déduire du théorème précédent
l'existence d'une version continue du processus
$(W_\rho^\hyp (z)\, :\ (z,\rho)\in \D\times ]0,+\infty[)$.
Pour une telle version continue, la famille des moyennes
du champ libre gaussien
sur les cercles du disque de Poincaré centrés en un point donné réalise un mouvement brownien réel standard de la façon suivante.
Pour tout $z_0$ dans $\D$,
le processus $(B_t (z_0)\, :\ t\in \R_+)$ défini par
$$
B_t (z_0)=
\left\{
\begin{array}{ll}
W^\hyp_{\rho (t)} (z_0) \quad & \mbox{si $t>0$}\\
0 & \mbox{si $t=0$}
\end{array}
\right.
\quad\quad
\mbox{où $\rho (t)={\rm argtanh} (e^{-t})$}
$$
est un mouvement brownien réel
standard
puisque l'on a, pour tous $s$ et $t$ dans $\R_+$,
$$
\cov \left(B_s (z_0) , B_t (z_0)\right)
=
\cov \left(W_{\rho (s)} (z_0),W_{\rho (t)} (z_0)\right)
=-\ln \max(e^{-s},e^{-t})
=\min(s,t)
\, .
$$

Notons que, compte tenu des limites:
$$
\lim_{t\to 0^+} \rho (t)=+\infty
\qeq
\lim_{t\to +\infty} \rho (t)=0\, ,
$$
les cercles de petits rayons donnent le mouvement brownien dans le temps long,
tandis que ceux de grands rayons correspondent au mouvement brownien dans un voisinage de l'origine.

\end{document}